\documentclass[12pt]{article}
\usepackage{amssymb,amsfonts,amsthm, amsmath}
\usepackage{color}

\newtheorem{thm}{Theorem}[section]
\newtheorem{cor}[thm]{Corollary}
\newtheorem{lem}[thm]{Lemma}
\newtheorem{prop}[thm]{Proposition}

\newtheorem{rem}[thm]{Remark}

\def\eop{\hfill\rule{2.5mm}{2.5mm}}
\def\pf{\par\smallbreak\noindent {\bf Proof.} \ }
\DeclareMathOperator*{\esssup}{ess\,sup}

\begin{document}

\title{
{\textbf{\Large{
Asymptotic analysis on a non-standard Hilbert space of non-absolutely integrable functions }}} \vspace{-4pt}
\author{
\textsc{F. Andrade da Silva}\,\footnote{Grant \#2021/12213-5, São Paulo Research Foundation (FAPESP).}
\,\,\&\,\,
\textsc{K. Gonzalez}\,\footnote{Grant: Morá Miriam Rozen Gerber, Weizmann Institute of Science.}
\,\,\&\,\,
\textsc{T. Jord\~{a}o}\,\footnote{Grant \#2022/11032-0, São Paulo Research Foundation (FAPESP).}
}}
\date{}

\maketitle 

\vspace{-30pt}
\begin{center}
\parbox{13 cm}{{\small {\bf Abstract:} 
In this work, we study the Kuelbs-Steadman-2 space (KS-2 space), a Hilbert space constructed via the Henstock-Kurzweil integral, which allows handling non-absolutely integrable functions. We present the construction of the KS-2 space over measurable subsets of $\mathbb{R}^d$ and explore its functional properties with particular focus on integral operators associated with symmetric kernels. A Mercer-type representation theorem is established for such kernels in a KS-2 space, leading to the characterization of the associated Reproducing Kernel Hilbert Spaces (RKHS). As an application, we derive asymptotic upper and lower bounds for the covering numbers of the embedding of the RKHS into the KS-2 space, highlighting how the Fourier coefficients decay rate of the kernels influences the estimates. }}

\end{center}

\thispagestyle{empty}

\section{Introduction}\label{intro}

In the mid-1950s, Ralph Henstock and Jaroslav Kurzweil came up with a new and simple formulation of integral which encompasses the Riemann, Newton, and Lebesgue integrals (see \cite{Henstock-First, Kurzweil-First}). In their definition, the intuitive approach of the Riemann integral is preserved, but unlike the latter, Henstock and Kurzweil considered a strictly positive function $\delta$ (called gauge) to calibrate the length of each subinterval of the domain. By making this ``small adjustment'', it turned out that their integral, known as the Henstock-Kurzweil integral (HK-integral, for short), recovers all primitives as integrals. Moreover, this integral allows us to deal with highly oscillating integrands. The KH-integral can be extended to functions defined in a complete metric measure space endowed with a radon measure (see \cite{bongiorno, gill2016functional} for instance). 

We denote by $\mathcal{P}=\left\{([t_{i-1}, t_i], \tau_i): i=1,\ldots, n \right\}$, with $t_{i-1}<t_i$ and $\tau_i\in [t_{i-1}, t_i]$ for all $i\in \{1,\ldots, n\}$, a tagged partition of the interval $[a, b] \subset \mathbb{R}$ ($a<b$). For a function $\delta:[a, b] \rightarrow(0, \infty)$, a \textit{HK-$\delta$ partition} is a tagged partition $\mathcal{P}=\left\{([t_{i-1}, t_i], \tau_i): i=1,\ldots, n \right\}$ calibrated by $\delta$, i.e., it holds that $0< t_i-t_{i-1}<\delta(\tau_i)$, for $1\leq i \leq n$. In this case, the function $\delta$ is called a \textit{gauge} on $[a,b]$. It is not hard to see that for any gauge $\delta$ the tagged partition $\mathcal{P}$ is a HK-$2\delta$ partition if and only if it holds that $[t_{i-1}, t_i]\subset (\tau_i-\delta(\tau_i), \tau_i+\delta(\tau_i))$, for $i=1,\ldots, n$. A function $f:[a, b] \rightarrow \mathbb{R}$ is \textit{HK-integrable} if there exists $\mathrm{HK}(f)\in\mathbb{R}$ such that for all $\epsilon>0$, there exists a gauge $\delta:[a, b] \longrightarrow(0, \infty)$ such that for all HK-$\delta$ partition $\mathcal{P}=\left\{([t_{i-1}, t_i], \tau_i): i=1,\ldots, n \right\}$, it holds that
\begin{equation}\label{def-HK-integral}
\left|\sum_{i=1}^n f\left(\tau_i\right)(t_i-t_{i-1}) -\mathrm{HK}(f)\right|<\epsilon.    
\end{equation}
In this case, $\mathrm{HK}(f)$ is called the Henstock-Kurzweil integral, or shortly, the \textit{HK-integral} of $f$, and we write
$$
\mathrm{HK}\int_a^b f(t) d t:=\mathrm{HK}(f).
$$
It is easy to see that if $f$ is integrable in the sense of Riemann, then $f$ is HK-integrable. In fact, for any $\epsilon>0$, if we consider $\gamma>0$ such that 
$$
\left|\sum_{i=1}^n f\left(\tau_i\right)(t_i-t_{i-1})-\int_{a}^bf(t)dt\right|<\epsilon,
$$
for all tagged partition $\mathcal{P}$ with $|\mathcal{P}|=\max\{t_i-t_{i-1}: i=1,\ldots, n\}<\gamma$, then the constant function $\delta: [a,b]\rightarrow(0,\infty)$ given by $\delta\equiv (b-a)/\eta$, with $\eta\in\mathbb{N}$ such that $(b-a)/\eta <\gamma$, is a gauge that fulfills the condition (\ref{def-HK-integral}).

Let $L([a,b]):=L([a,b], \mathcal{L}, \lambda)$ the vector space of Lebesgue integrable functions on $[a,b]$. We write $\int_{a}^bf(t)dt$ for the Lebesgue integral of $f\in L([a,b])$. The relation of the HK-integral with the sense of integral of Lebesgue is established bellow. The reader can find the proofs (omitted here) in \cite[Theorem 3.14]{gill2016functional} and in \cite[Theorem~4.60 and Corollary 4.62]{theories}. Consider $f$ an HK-integrable function on $[a,b]$, then the following properties hold true:\\
i) if $f$ is a bounded or a non-negative function, then $f\in L([a,b])$,\\
ii) $|f|$ is HK-integrable if, and only if, $F(t)=\int_a^t f(s)ds$ is of bounded variation on $[a,b]$, and \\  
iii) if $g$ is  HK-integrable on $[a,b]$ and $|f|\leq g$, then $f\in L([a,b])$.

Due to the fact that conditions over the integrands can be weakened by considering more general integrals (see \cite{reid1975anatomy}), the space of generalized ordinary differential equations (we write generalized ODEs for short) provides a framework to deal with many differential and integral equations involving non-absolute integrable functions, once solutions of generalized ODEs are described by Kurzweil integral. Hence, the theory of generalized ODEs has been shown to act as a unifying theory for many equations. For instance, measure functional differential equations, equations with impulses, dynamic equations on time scales, and integral equations can be regarded as generalized ODEs (see, e.g., \cite{B-F-M} and the references therein). Results on stability through Lyapunov functionals are presented in \cite{converse-regular, converse-uniform}. The study of generalized ODE has shown great relevance in the field of differential equations, and its investigation has increased in recent years. In fact, generalized ODEs became a theme at the Mathematics Subject Classification (MSC2020), namely, 34A06. 

In this work, we present the concept of the HK-integral on the $d$-dimensional Euclidean space and for Lebesgue mensurable subsets of $\mathbb{R}^d$ in order to introduce the Kuelbs-Steadman-2 space (KS-2 space). The KS-2 space for $d=3$, denoted by $KS^2(\mathbb{R}^3)$ is an important space used by the authors in \cite{gill2016functional} to provide a rigorous foundation for the Feynman formulation of quantum mechanics. We briefly explore the theory of integral operators generated by kernels in the Kuelbs-Steadman-2 space to establish a Mercer series representation for symmetric kernels and to present the reproducing kernel Hilbert space related to the definite positive kernels. As an application of the theory, we present an upper bound for the covering numbers of the unit ball for certain classes of Reproducing Kernel Hilbert Space (RKHS) of kernels in the Kuelbs-Steadman-2 space. The technique is based on the application of standard properties for covering numbers, and it was successfully applied in \cite{GONZALEZ2024128121} to obtain upper and lower estimates for the covering numbers of spherical continuous zonal, and positive definite kernels. 

The covering number or Kolmogorov $\epsilon$-entropy \cite{kolmogorov-MR0112032} is a popular concept that plays an important role in applied areas such as kernel-based learning algorithms and Gaussian process \cite{Li-MR1733160, Minh2010-MR2677883, WSS-MR1873936}, which are fundamental techniques in computer and data sciences.  The bounds for the covering numbers are a useful tool to estimate the probabilistic error of the statistical nature of the observations from which the algorithms designed are learning from \cite{smale-MR1864085, Zhou-2003-MR1985575}. Furthermore, building on recent advances in the theory of non-absolutely integrable functions, our framework based on the HK-integral and KS-2 space has the potential to impact a wide range of scientific and engineering research fields. These include the physical sciences (e.g., quantum mechanics via Feynman path integrals), biological modeling (e.g., impulsive differential equations), , and engineering (e.g., control theory for systems exhibiting discontinuous behavior). 

This paper is organized as follows. In Section \ref{integrability}, discusses the relationship between the Lebesgue and the HK notions of integrability. Subsection \ref{d-dimensional-integral} introduces the concept of the $d$-dimensional HK-integral defined on Lebesgue measurable subsets of $\mathbb{R}^d$. Section \ref{sec-KS-space} presents the KS-2 space of functions on $\mathbb{R}^d$ and establishes an orthonormal complete system for the Hilbert space $\mathrm{KS}^{\infty}(\mathbb{R}^d)$ with inner product $ \langle\,\cdot , \cdot \,\rangle_{KS2}$. In Section \ref{Mercer-representation}, we explore Mercer’s theorem for symmetric kernels in $\mathrm{KS}^2(\mathbb{R}^d)$ and provide a representation theorem for the associated reproducing kernel Hilbert space. Subsection \ref{covering} establishes an upper bound for the covering numbers of the embedding $\mathcal{H}_K \hookrightarrow \mathrm{KS}^2(\mathbb{R}^d)$. Using the Mercer representation of the kernel $K$ with exponentially decaying eigenvalues, it derives an explicit asymptotic upper bound for the growth rate of the covering numbers as $\epsilon \rightarrow 0^+$, providing a quantification of the compactness of the embedding. Subsection \ref{lower} provides a matching lower bound for the same embedding, showing that 
\[
\log(\mathcal{C}(\epsilon,I_K ))\asymp {\log\left(\frac{1}{\epsilon}\right)}^{2} .
\]
as $\epsilon \rightarrow 0^+$. It describes the asymptotic behavior for the covering numbers in the setting of this work.

\section{Integrability: $\mathrm{HK}$ versus Lebesgue}
\label{integrability}

We start this section by presenting a result that reflects the relation between the Lebesgue and the HK senses of integral. Furthermore, it ensures us that HK-integral has a weaker sense of integrability since we deliver a nice and simple example of a function that is HK-integrable and that is not Lebesgue integrable.

\begin{prop}\label{HK-integral}
If $f\in L([a,b])$, then $f$ is $\mathrm{HK}$-integrable and it holds that
$$
\mathrm{HK}\int_a^b f(t) dt=\int_a^b f(t) dt.
$$
The reciprocal does not hold true. Furthermore, a $\mathrm{HK}$-integrable function cannot be written necessarily as a sum of a Lebesgue with an improper Riemann integrable function.
\end{prop}

\pf The first part of the proof is an standard result of the HK-integral theory and it was omitted. We present a simple example of  a HK-integrable function on $[0,1]$ that is not in $L([a,b])$. The example is given in terms of an improper Riemann integrable function as follows. Consider $h:[0,1]\to \mathbb{R}$ defined by $h(0)=0$, and 
$$
h(t) = 2t\cos\left(\pi/t^2\right)+ \frac{2\pi}{t} \sin\left(\pi/t^2\right), \quad  t\in (0,1].
$$
Consider now the Dirichlet function $D_h: [0,1]\to \mathbb{R}$ associated to $h$, given by 
$$
D_h(t) = \left\{ \begin{array}{ll}
h(t) +1, & t\in [0,1]\cap \mathbb{Q},\\
h(t), & t\in (0,1]\cap(\mathbb{R  \setminus \mathbb{Q}}).
\end{array}
\right.
$$
Clearly, $D_h$ is not improper Riemann integrable or integrable in the sense of Lebesgue but $D_h$ is a HK-integrable function. 

Finally, we consider $f:[0,1]\to \mathbb{R}$ defined by $f(0)=0$, and
$$
f(2^{-n}(1+x)) = (-1)^{n+1}\frac{2^n}{n}, \quad  x\in (0,1],
$$ 
for $n=1,2,\ldots$. It is clear that $f$ is well-defined. For $n=1,2,\ldots$, we have that 
\begin{eqnarray*}
\lim_{c\to 0^+} \int_c^1 f(2^{-n}(1+x)) dx =  \lim_{c\to 0^+}  (-1)^{n+1}\frac{2^n}{n} \int_c^1 dx &=&   \lim_{c\to 0^+} (-1)^{n+1}\frac{2^n}{n} (1-c) \\ &=&  (-1)^{n+1}\frac{2^n}{n}.
\end{eqnarray*}
Also, by the Hake's Theorem (see e.g. \cite[Theorem~12.8]{Modern}) $f(2^{-n}(1+x))$ is HK-integrable on $[0,1]$ and its HK-integral is 
\begin{equation}\label{1}
 \mathrm{HK}\int_0^1 f(2^{-n}(1+x))  dx  =  (-1)^{n+1}\frac{2^n}{n} .
\end{equation}
By the Substitution Theorem (see \cite[Theorem~4.60]{Modern}, for instance), we get 
\begin{equation}\label{2}
\mathrm{HK}\int_0^1  f(2^{-n}(1+x)) 2^{-n} dx = \mathrm{HK}\int_{2^{-n}}^{2^{-n+1}} f(x)dx. 
\end{equation}
By equations (\ref{1}) and (\ref{2}), we deduce that
$$
\mathrm{HK}\int_0^1  f(x)dx  =\sum_{n=1}^{\infty } \mathrm{HK}\int_{2^{-n}}^{2^{-n+1}} f(x)dx =  \sum_{n=1}^{\infty } \frac{(-1)^{n+1}}{n} =\ln(2).
$$
Similarly, $g:[0,1]\to \mathbb{R}$ given by $g(2^{-n}(1+x))= f(x)$ for $x\in (0,1]$ and $g(0)=0$ is Henstock-Kurzweil  integrable on $[0,1]$. If we assume that $g=l+k$, where $l$ is a Lebesgue integrable function and $k$ is an improper Riemann integrable function, then there exists a partition $P=\{ t_j: t_{j-1}< t_j, \quad j=1,\ldots, m\}$ of $[0,1]$ such that $k$ is an improper Riemann integrable function on $[x,y]$ for all $0=t_0 < x< y < t_{1}$. Let $n$ be large enough so that $[2^{-n}, 2^{-n+1}]\subset  (0, t_1)$. Since $|g| \leq |l| + |k|$ and $|l| + |k|$ is  Henstock-Kurzweil  integrable on  $[2^{-n}, 2^{-n+1}]$, it follows from item (iv) that $g$ is Lebesgue integrable on $[2^{-n}, 2^{-n+1}]$ and, consequently, $f$ is Lebesgue integrable. However, 
\[ \int_0^1 |f(x)|dx \geq \sum_{n=1}^N \frac{1}{n} \to \infty\]
which is a contradiction. This shows that $g$ is not a sum of Lebesgue integrable and improper Riemann integrable function. 
\eop

The Henstock-Kurzweil integral can be defined on mensurable sets of $\mathbb{R}$ as follows. For $E\in\mathcal{L}$, a family $\mathcal{P}=\left\{(J_i, \tau_i): i=1,\ldots, n \right\}$ of closed intervals $J_i=[t_{i-1}, t_i]$ and points $\tau_i\in J_i\cap E$, $i=1,\ldots, n$, satisfying
$\lambda(J_i \cap J_j)=0$, for any $i \neq j$, and $E=\cup_{i=1}^m (J_i\cap E)$ is called a \textit{tagged partition of $E$}. For a function $\delta: E\rightarrow (0,\infty)$, a tagged partition $\mathcal{P}=\left\{(J_i, \tau_i): i=1,\ldots, n \right\}$ of $E$ is a HK-$\delta$ partition if $[t_{i-1}, t_i]\subset (\tau_i-\delta(\tau_i), \tau_i+\delta(\tau_i)), \quad i=1,\ldots, n$. For $E\in\mathcal{L}$ the function $f:E \rightarrow \mathbb{R}$ is HK-integrable if there exists $\mathrm{HK}(f)\in\mathbb{R}$ such that for all $\epsilon>0$, there exists a gaude $\delta:E \longrightarrow(0, \infty)$ such that for all HK-$\delta$ partition $\mathcal{P}$ of $E$, with $[t_{i-1}, t_i]\subset (\tau_i-\delta(\tau_i), \tau_i+\delta(\tau_i))$, $i=1,\ldots, n$, it holds that
\begin{equation}\label{def-HK-integral-mensurable}
\left|\sum_{i=1}^n (t_i-t_{i-1}) f\left(\tau_i\right)-\mathrm{HK}(f)\right|<\epsilon.    
\end{equation}
In this case, $\mathrm{HK}(f)$ is called the HK-integral of $f$, and we write
$$
\mathrm{HK}\int_E f(t) d t:=\mathrm{HK}(f).
$$
We highlight that if $f$ is HK-integrable on any Lebesgue mensurable set $E\in\mathcal{L}$, then $f$ is Lebesgue integrable on $E$ (see \cite[Theorem 3.14]{gill2016functional}).

\subsection{The $d$-dimensional HK-integral}
\label{d-dimensional-integral}

In this subsection, we present the concept of the $d$-dimensional HK-integral defined on Lebesgue-mensurable subsets of $\mathbb{R}^d$. The presentation of the concept of the HK-integral will be brief, and for details, we suggest the source \cite[p. 112-123]{gill2016functional} and references cited in there.

For $E\in\mathcal{L}$ a measurable subset of $\mathbb{R}^d$, we keep the notation $L(E):=L(E, \mathcal{L}, \lambda)$ for the vector spaces of Lebesgue integrable functions on $E$, and we employ the notation $\int_Ef(x)dx$ for the Lebesgue integral of $f\in L(E)$. We consider $\mathbb{R}^d$, for $d\geq 1$, with the maximum norm 
$$
\|x\|=\max\{|x_i| : 1 \leq i \leq d\}, \quad x=(x_1, \ldots, x_d)\in\mathbb{R}^d.
$$ 
The closed ball centered at $x\in\mathbb{R}^d$ with radius $r>0$, denoted by $D(x, r)$, here is a cube as follows 
$$ 
D(x, r)=\prod_{i=1}^{d}\left[a_i, b_i\right],\quad \mbox{with} \,\, \max\{|b_i-a_i|: 1\leq i\leq d\}= 2r.
$$ 
Similarly, the open ball center at $x\in\mathbb{R}^d$ with radius $r>0$ is given by
$ 
B(x, r)=\prod_{i=1}^{d}\left(a_i, b_i\right)$, with $\max\{|b_i-a_i|: 1\leq i\leq d\}< 2r$. 
Moreover, we write 
$(J, x) = \prod_{i=1}^{d}\left[a_i, b_i\right]$, with $x_i\in [a_i, b_i]$, for $1\leq i\leq d$, and we refer to $J$ as closed interval of $\mathbb{R}^d$. 

A family $\mathcal{P}=\left\{(J_j, x_j): j=1,\ldots, m \right\}$ of closed intervals $J_j$ of $\mathbb{R}^d$ and points $x_j\in J_j\cap E$ , $j=1,\ldots, m$, satisfying
$\lambda(J_j \cap J_k)=0$, for any $j \neq k$, and $E=\cup_{i=1}^m (J_i\cap E)$ is called a \textit{tagged partition of $E$}. For $\delta:E \longrightarrow(0, \infty)$ a function defined on the measurable set $E\subset \mathbb{R}^d$, we say that a tagged partition of $\mathcal{P}=\left\{\left(J_j, x_j\right): x_j \in J_j, \,\, 1 \leq j \leq m\right\}$ of $E$ is a \textit{HK-$\delta$ partition} for $E$ if,
\[
J_j\subset B(x_j, \delta(x_j)), \quad \mbox{for all}\,\,  j=1, \ldots, m.
\]
In this case, the function $\delta$ is also called a gauge on $E$, and for every compact set $C \subset \mathbb{R}^d$, one can always to find a HK-$\delta$ partition for $C$ (see \cite[Lemma 3.21]{gill2016functional}). 

For $E\in\mathcal{L}$, a function $f: E \rightarrow \mathbb{R}$ is \textit{HK-integrable} on $E$, if there exists $\mathrm{HK}_E(f)\in\mathbb{R}$ such that for any $\epsilon>0$, there is a gauge $\delta$ on $E$ and a HK-$\delta$ partition of $E$ such that
$$
\left|\sum_{i=1}^m f\left(x_i\right) \lambda(J_i)-\mathrm{HK}_E(f)\right|<\epsilon.
$$
If $f: E \rightarrow \mathbb{R}$ is HK-integrable, the number $\mathrm{HK}_E(f)$ is the HK-integral of $f$ on $E$, and we write
$$
\mathrm{HK}\int_E f(x) dx:=HK_E(f).
$$

If $f: E \rightarrow \mathbb{R}$ is such that $f\in L(E)$, then $f$ is HK-integrable on $E$ (see \cite[Theorem 3.24]{gill2016functional}). Additionally, if $f\in L(E)$, then for any $\varepsilon>0$ and any $x\in E$, there is an open set $G(x)$ containing $x$ such that for all family of closed sets contained in $E$, $\left\{B_1, B_2, \ldots\right\}$ satisfying $\lambda(B_i \cap B_j)=0$, for $i \neq j$, $E=\bigcup_{i=1}^{\infty} B_i$, $\lambda\left(E \backslash \cup_{k=1}^{\infty} B_k\right)=0$, and $\{x_1, x_2, \ldots\}$ fulfilling $x_k \in B_k \subset G(x_k)$, for $k=1,2,\ldots$, we have
$$
\left|\sum_{k=1}^{\infty} f(x_k) \lambda(B_k)-\int_E f(x) dx\right|<\varepsilon.
$$
In this case, clearly, it holds that $\int_E f(x) dx=\mathrm{HK}\int_E f(x) dx$.

We write $\mathrm{HK}(\mathbb{R}^d)$ for the vector space of HK-integrable functions on $\mathbb{R}^d$ endowed with the norm 
$$
\|f\|_{\mathrm{HK}}=\sup _{r>0}\left|\mathrm{HK}\int_{D(0,r)} f(x) dx\right|<\infty.
$$
The HK-integral is equivalent to the restricted Denjoy integral, thus the formula above defines a norm for both $D(\mathbb{R}^d)$ and the wide sense Denjoy integrable functions (see e.g. \cite[Chapter~11]{Gordon}).

\section{The Kuelbs-Steadman Space (KS-2 space)}\label{sec-KS-space}

In this section, we present the KS-2 space of functions on $\mathbb{R}^d$, the main reference here is  \cite{gill2016functional}. For $1\leq p\leq\infty$ we write $(L^p(\mathbb{R}^d), \|\,\cdot \,\|_p)$ for the Lebesgue spaces of all real-valued functions $f$ on $\mathbb{R}^d$ such that $|f|^p\in L(\mathbb{R}^d)$, for $p\neq \infty$, and for $p=\infty$ it is the space of all functions $f$ on $\mathbb{R}^d$ such that $\|f\|_{\infty}=\esssup\{|f(x)|: x\in\mathbb{R}^d\}<\infty$. For $E\in\mathcal{L}$ a measurable subset of $\mathbb{R}^d$, we keep the notation $L(E):=L(E, \mathcal{L}, \lambda)$ for the vector spaces of Lebesgue integrable functions on $E$, and we employ the notation $\int_E f(x)dx$ for the Lebesgue integral of $f\in L(E)$. 

Let $\mathbb{Q}^d\subset\mathbb{R}^d$ arranged as $\mathbb{Q}^d=\left\{q_1, q_2, \cdots\right\}$. Let $B_n(q_i)$ be the closed ball centered at $q_i$ with radius $1/2^{n}\sqrt{d}$, for $i=1,2,\ldots$. The edges of the cubes $B_n(q_i)$ are parallel to the coordinate axes and length $e_n=2/2^n\sqrt{d}$, for $n=1,2,\ldots$. Organizing the cubes $\{B_n(q_i)\}_{i,n=1}^{\infty}$ by considering the lexicography order on the pairs $(n,i)$, for $n,i=1,2\ldots$, we consider the resultant sequence written as $\left\{B_k : k=1,2,\ldots\right\}$. If the notation $\chi_{B_k}$ stands for the characteristic function of $B_k$ on $\mathbb{R}^d$, for $k=1,2\ldots$, then it is clear that $\chi_{B_k}\in L^p(\mathbb{R}^d) \cap L^{\infty}(\mathbb{R}^d)$, since
\begin{equation}\label{volume-bola}
\int_{\mathbb{R}^d}\chi_{B_k}  d x =\lambda(B_k)=\frac{\pi^{d/2}}{2^{d-1}d^{d/2+1}\Gamma(d/2)}\frac{1}{2^{dk}}\leq \frac{\pi^{d/2}}{2^{d-1}d^{d/2+1}\Gamma(d/2)},\qquad k=1,2\ldots. 
    \end{equation}
The symbol $\Gamma$ stands for the Gamma function, the Stirling's approximation formula ($\Gamma(x)\approx (2\pi/x)^{1/2}(x/e)^x$ as $x\to\infty$ and for $x>0$) implies  
$$
\Gamma(d/2)\approx \frac{(2\pi)^{1/2}(d/2)^{d/2}}{(d/2)^{1/2}e^{d/2}}=\left(\frac{2\pi}{e^d}\right)^{1/2}\left(\frac{d}{2}\right)^{(d-1)/2}, \quad \mbox{as $d\to\infty$}.
$$
Above the notation $a_n\approx b_n$ stands for $a_n/b_n\to 1$, as $n\to\infty$. Then, there exists a constant $c$ does not depending on the dimension $d$ such that
$$
\int_{\mathbb{R}^d}\chi_{B_k}  d x =\lambda(B_k)\leq \frac{c\pi^{d/2}}{2^{d-1}d^{d/2+1}}\left(\frac{e^d}{2\pi}\right)^{1/2}\left(\frac{2}{d}\right)^{(d-1)/2}=\frac{ce^{d/2}}{\pi d^{d-1}}\left(\frac{\pi}{2}\right)^{d/2}.
$$

For $k=1,2,\ldots$ and $1\leq p\leq\infty$, define $F_k$ on $L^p(\mathbb{R}^d)$ by
\begin{equation}
F_k(f)=\int_{\mathbb{R}^d} \chi_{B_k}(x) f(x) dx, \qquad f\in L^p(\mathbb{R}^d), \quad \mbox{if $p\neq\infty$},    
\label{Fk-LP}
\end{equation}
and
\begin{equation} 
F_{k}(f)=\esssup\left\{\chi_{B_k}(x) f(x): x\in \mathbb{R}^d\right\}, \qquad f\in L^{\infty}(\mathbb{R}^d), \quad \mbox{if $p=\infty$}.
\label{Fk-Linfty}
\end{equation}
The functional linear $F_k$, defined in (\ref{Fk-LP}) or (\ref{Fk-Linfty}), is bounded on $L^p(\mathbb{R}^d)$, and $\|F_k\| \leq 1$, for $k=1,2\ldots$. Furthermore, if $f\in L^p(\mathbb{R}^d)$ and $F_k(f)=0$ for $k=1,2,\ldots$, then  $f=0$. 

We consider $\mu$ the measure defined for all the mensurable subsets of $\mathbb{R}^d \times \mathbb{R}^d$ (endowed with the Lebesgue product measure) given by
\begin{equation}\label{def-medidamu}
d\mu(x,y)=\left(\sum_{k=1}^{\infty} 2^{-k} \chi_{B_k}(x) \chi_{B_k}(y)\right) dx dy .    
\end{equation}
It is not hard to check that the volume $V_{\mu}$ of $\mathbb{R}^d\times\mathbb{R}^d$ under the measure $\mu$ is finite. In fact,
\begin{eqnarray*}
V_{\mu}=\sum_{k=1}^{\infty} 2^{-k} \int_{\mathbb{R}^d\times\mathbb{R}^d}\chi_{B_k}(x) \chi_{B_k}(y) dx dy =\sum_{k=1}^{\infty} 2^{-k} (\lambda(B_k))^2=\frac{\alpha^2_d}{2^{2d+1}-1},
    \end{eqnarray*}
with $\lambda(B_k)= \alpha_d2^{-kd}$, for $k=1,2\ldots$, and $\alpha_d:= \pi^{d/2}/(2^{d-1}d^{d/2+1}\Gamma(d/2))$.

We consider the product $\langle\, \cdot, \cdot\,\rangle_{KS2}: L^1(\mathbb{R}^d)\times L^1(\mathbb{R}^d)\rightarrow\mathbb{R}$ by
\begin{equation*}
\langle f, g\rangle_{KS2}  =\int_{\mathbb{R}^d \times \mathbb{R}^d} f(x) g(y) d \mu(x,y), \qquad f,g\in L^1(\mathbb{R}^d).
\end{equation*}
The product $\langle\, \cdot, \cdot\,\rangle_{KS2}$ is an inner product on $L^1(\mathbb{R}^d)$ and it is given by \begin{equation}\label{def-innermu}
\langle f, g\rangle_{KS2}  = \sum_{k=1}^{\infty} 2^{-k}F_k(f)F_k(g), \qquad f,g\in L^1(\mathbb{R}^d).
\end{equation}
The completion of $(L^1(\mathbb{R}^d), \langle\,\cdot ,\cdot \,\rangle_{KS2})$, denoted by $(\mathrm{KS}^2(\mathbb{R}^d), \langle\,\cdot , \cdot \,\rangle_{KS2})$, is called \textit{KS-2 space} (Kuelbs-Steadman-2 space). The completion of $(L^{\infty}(\mathbb{R}^d), \langle\,\cdot ,\cdot \,\rangle_{KS2})$, denoted by $(\mathrm{KS}^{\infty}(\mathbb{R}^d), \langle\,\cdot , \cdot \,\rangle_{KS2})$, is called KS-$\infty$-space (Kuelbs-Steadman-$\infty$ space).

\begin{lem}\label{systemofKS2}
The family of normalized characteristic functions 
$$
\beta=\left\{Y_k:=\frac{2^{(k-1)/2}}{\lambda(B_k)}\chi_{B_k} : k=1,2,\ldots \right\}
$$
is an orthonomal complete system of the Hilbert space $(\mathrm{KS}^2(\mathbb{R}^d), \langle\,\cdot \,, \,\cdot \,\rangle_{KS2})$.    
\end{lem} 

The Lemma \ref{systemofKS2} and an application of the Parseval's identity \cite[p.175]{FollandReal} ensures us 
\begin{equation}\label{FseriesKS2}
f(x)=\sum_{k=1}^{\infty}\langle f, Y_k\rangle_{KS2}Y_k(x), \qquad  x\in\mathbb{R}^d, \mu-a.e.,
\end{equation}
and
\begin{equation}\label{FseriesKS2-norm}
\|f\|^2_{KS2}=\sum_{k=1}^{\infty}|\langle f, Y_k\rangle_{KS2}|^2, \qquad f\in \mathrm{KS}^2(\mathbb{R}^d).
\end{equation}

The Lebesgue space $L^1(\mathbb{R}^d)$ is densely and continuously embedded in $\mathrm{HK}(\mathbb{R}^d)$, defined in the previous section. This makes possible to embed $L^1(\mathbb{R}^d)$ densely and continuously in the Hilbert space $\mathrm{KS}^2(\mathbb{R}^d)$ in the following way. For any $f \in \mathrm{HK}(\mathbb{R}^d)$, it holds
$$
\|f\|_{K S^2}^2 =\sum_{k=1}^{\infty} 2^{-k}\left|\int_{\mathbb{R}^d} \chi_{B_k}(x) f(x) dx\right|^2
 \leqslant \sup _{k\geq 1}\left|\int_{\mathbb{R}^d} \chi_{B_k}(x) f(x) dx\right|^2 \leqslant\|f\|_{\mathrm{HK}}^2.
$$
In symbols, $L^1(\mathbb{R}^d)\hookrightarrow \mathrm{HK}(\mathbb{R}^d)\hookrightarrow \mathrm{KS}^2(\mathbb{R}^d)$. It is also known that $L^2(\mathbb{R}^d)$ is a dense subspace of $\mathrm{KS}^2(\mathbb{R}^d)$ \cite[Theorem 3.25]{gill2016functional}, and the following properties hold. If $f_n \rightarrow f$ weakly in $L^2(\mathbb{R}^d)$, then $f_n \rightarrow f$ strongly in $KS^2(\mathbb{R}^d)$, $KS^2(\mathbb{R}^d)$ is uniformly convex and reflexive, since the dual space of $KS^2(\mathbb{R}^d)$ is $K S^2(\mathbb{R}^d)$, and $K S^{\infty}(\mathbb{R}^d) \subset KS^2(\mathbb{R}^d)$ \cite[Theorem 3.27]{gill2016functional}.

\subsection{Integral operators on the KS-2 space} 
\label{integral-operators}

We consider on $\mathbb{R}^d\times \mathbb{R}^d$ the product Lebesgue measure ($dz=dxdy$, $z=(x,y)\in \mathbb{R}^d\times \mathbb{R}^d$) to define the standard Lebesgue space $L^p(\mathbb{R}^d\times\mathbb{R}^d)$, for $1\leq p\leq\infty$. For $k,j=1,2,\ldots$, let $\chi_k\otimes \chi_j: \mathbb{R}^d\times \mathbb{R}^d\rightarrow\mathbb{R}$ be the function given by 
$$
(\chi_k\otimes \chi_j)(z)=\chi_k(x)\chi_j(y), \qquad z=(x,y)\in\mathbb{R}^d\times\mathbb{R}^d.
$$ 
Define $F_{k,j}$ on $L^p(\mathbb{R}^d\times\mathbb{R}^d)$ by
\begin{equation*}
F_{k,j}(f)=\int_{\mathbb{R}^d\times\mathbb{R}^d} (\chi_{B_k}\otimes\chi_{B_j})(z) f(z) dz, \qquad f\in L^p(\mathbb{R}^d\times \mathbb{R}^d), \quad \mbox{if $p\neq\infty$},    
\end{equation*}
and, if $p=\infty$, then 
\begin{equation*} 
F_{k,j}(f):=\esssup\left\{(\chi_{B_k}\otimes\chi_{B_k})(z) f(z): z=(x,y)\in \mathbb{R}^d\times\mathbb{R}^d\right\}, \qquad f\in L^{\infty}(\mathbb{R}^d\times \mathbb{R}^d).
\end{equation*}
It is not hard to see that the linear functional $F_{k,j}$ is bounded on $L^p(\mathbb{R}^d\times \mathbb{R}^d)$, satisfying $\|F_{k,j}\| \leq 1$, for $k,j=1,2\ldots$, and if $F_{k,j}(f)=0$, for $k,j=1,2,\ldots$, then  $f=0$, for any $f\in L^p(\mathbb{R}^d\times\mathbb{R}^d)$. Analogously to the previous case, endow $(\mathbb{R}^d\times \mathbb{R}^d)^2$ with the product measure $\mu\times\mu$, satisfying
\begin{equation*}
d(\mu\times\mu)(z,w)=\left(\sum_{k,j=1}^{\infty} 2^{-k-j} (\chi_{B_k}\otimes\chi_{B_k})(z)(\chi_{B_j}\otimes\chi_{B_j})(w)\right) dzdw,    
\end{equation*}
for $z=(s,t),w=(x,y)\in(\mathbb{R}^d\times \mathbb{R}^d)^2$. Consider the completion of $L^1(\mathbb{R}^d\times \mathbb{R}^d)$ with the inner product 
\begin{equation*}
\langle f, g\rangle_{KS2}  =\int_{(\mathbb{R}^d \times \mathbb{R}^d)^2} f(z) g(w) d(\mu\times\mu)(z,w)  = \sum_{k,j=1}^{\infty} 2^{-k-j}F_{k,j}(f)F_{k,j}(g).
\end{equation*}
This completion delivers us the Hilbert space $(\mathrm{KS}^2(\mathbb{R}^d\times\mathbb{R}^d), \langle\,\cdot \,, \,\cdot \,\rangle_{KS2})$. 

Since $\beta=\{Y_k\}_k$ is an orthonomal complete system of the Hilbert Space $(\mathrm{KS}^2(\mathbb{R}^d), \langle\,\cdot \,, \,\cdot \,\rangle_{KS2})$ (Lemma \ref{systemofKS2}) we derive the following result.

\begin{lem}\label{systemofKS2-product}\cite[Theorem 65]{gill2009}
The family of $\otimes$-products of normalized characteristic functions 
$$
\gamma=\left\{Y_{k,j}:=Y_k\otimes Y_j : k, j=1,2,\ldots \right\},
$$
is an orthonomal complete system of the Hilbert Space $(\mathrm{KS}^2(\mathbb{R}^d\times\mathbb{R}^d), \langle\,\cdot \,, \,\cdot \,\rangle_{KS2})$.    
\end{lem}

Consider $K:\mathbb{R}^d\times \mathbb{R}^d\rightarrow\mathbb{R}$ a function of $\mathrm{KS}^2(\mathbb{R}^d\times\mathbb{R}^d)$. By Lemma \ref{systemofKS2-product}, the function $K$ can be written as
\begin{equation}\label{kernel}
K(x,y)=\sum_{k,j=1}^{\infty}\langle K, Y_{k,j}\rangle_{KS2}, Y_{k,j}(x,y), \qquad (x,y)\in\mathbb{R}^d\times \mathbb{R}^d, (\mu\times\mu)-a.e..
\end{equation}
Since $K(x, \,\cdot\, )\in \mathrm{KS}^2(\mathbb{R}^d)$, for any $x\in\mathbb{R}^d$, it is well defined the integral map
\begin{equation}\label{integral-operator}
I_K(f)(x):=\int_{\mathbb{R}^d\times\mathbb{R}^d}K(x,y)f(z)d\mu(y,z), \qquad f\in \mathrm{KS}^2(\mathbb{R}^d).  \end{equation}
By (\ref{def-innermu}), it is not hard to see that
$$
I_K(f)(x)=\langle K(x, \cdot)\,\, , f\rangle_{KS2}  = \sum_{k=1}^{\infty} 2^{-k}F_k(K(x,\, \cdot\, ))F_k(f), \qquad x\in \mathbb{R}^d,
$$
for any $f\in \mathrm{KS}^2(\mathbb{R}^d)$. Also, for any $x\in\mathbb{R}^d$, the Equation (\ref{kernel}) leads us to 
$$
\langle K(x,\,\cdot\,), Y_i\rangle_{KS2}=\sum_{k=1}^{\infty}\langle K, Y_{k,i}\rangle_{KS2}Y_k(x)=\sum_{k=1}^{\infty}2^{-k-i}F_{k,i}(K)Y_k(x), \qquad i=1,2,\ldots,
$$
due the fact that
$$
\langle K, Y_{k,i}\rangle_{KS2}=\sum_{l,j=1}^{\infty}2^{-l-j}F_{l,j}(K)F_{l,j}(Y_{k,i})=2^{-k-i}F_{k,i}(K), \qquad k,i=1,2,\dots. 
$$
Thus, for each $x\in\mathbb{R}^d$, it holds that
$$
\langle K(x,\,\cdot\,), Y_i\rangle_{KS2}=\sum_{k=1}^{\infty}\langle K, Y_{k,i}\rangle_{KS2}Y_k(x), \qquad i=1,2,\ldots,
$$
and, therefore,
\begin{eqnarray}\label{integral-operator-series}
I_K(f)(x)=\langle K(x, \cdot)\,\, , f\rangle_{KS2}=\sum_{j=1}^{\infty}\langle f, Y_j\rangle_{KS2}\sum_{k=1}^{\infty}\langle K, Y_{k,j}\rangle_{KS2}Y_k(x), \quad f\in \mathrm{KS}^2(\mathbb{R}^d).
\end{eqnarray}

\begin{rem}\label{remark-coef}
    Note that in the representation above the coefficients of $I_K(f)$, for $f\in \mathrm{KS}^2(\mathbb{R}^d)$, in terms of the orthonormal complete system $\{Y_k\}_k$ of $\mathrm{KS}^2(\mathbb{R}^d)$ are given by
$$
c_k= \sum_{j=1}^{\infty}\langle K, Y_{k,j}\rangle_{KS2}\langle f, Y_j\rangle_{KS2}, \qquad k=1,2,\ldots. 
$$
For $k=1,2,\ldots$, an application of the Cauchy-Schwarz for elements of $\ell_2$ implies 
\begin{eqnarray*}
|c_k|\leq  \left(\sum_{j=1}^{\infty}\left|\langle K, Y_{k,j}\rangle_{KS2}\right|^2\right)^{1/2}\|f\|_{KS2}.    
\end{eqnarray*}
Then, 
\begin{eqnarray*}
\sum_{k=1}^{\infty}|c_k|^2 \leq  \sum_{k=1}^{\infty}\sum_{j=1}^{\infty}\left|\langle K, Y_{k,j}\rangle_{KS2}\right|^2\|f\|_{KS2}^2= \|K\|_{KS2}^2\|f\|_{KS2}^2.   
\end{eqnarray*}
\end{rem}

The previous considerations and the Remark \ref{remark-coef} are summarized in the following result.

 \begin{prop}\label{prop-integral-operator} Let $K:\mathbb{R}^d\times \mathbb{R}^d\rightarrow\mathbb{R}$ be a function of $\mathrm{KS}^2(\mathbb{R}^d\times\mathbb{R}^d)$. The integral operator $I_K:\mathrm{KS}^2(\mathbb{R}^d)\rightarrow \mathrm{KS}^2(\mathbb{R}^d)$ given by 
$$
I_K(f)(x)=\langle K(x, \cdot)\,\, , f\rangle_{KS2}, \qquad x\in\mathbb{R}^d, \quad f\in \mathrm{KS}^2(\mathbb{R}^d), 
$$
is bounded and its operator norm satisfies $\|I_K\|\leq \|K\|_{KS2}$.
\end{prop}
    
The Proposition \ref{prop-integral-operator} could be also derived by density arguments based on the well-known property for integral maps defined on $L^1(\mathbb{R}^d)$ and for $K\in L^1(\mathbb{R}^d\times \mathbb{R}^d)$ (see \cite[Theorem 6.18]{FollandReal}, for example) and from the dense and continuous embedding $L^1(\mathbb{R}^d)\hookrightarrow \mathrm{HK}(\mathbb{R}^d)\hookrightarrow \mathrm{KS}^2(\mathbb{R}^d)$. These arguments appear disguised in a non-abstract way by the Parseval's identity we have applied a few times in the proof of Proposition \ref{prop-integral-operator}.  

Let $K:\mathbb{R}^d\times \mathbb{R}^d\rightarrow\mathbb{R}$ be a function of $\mathrm{KS}^2(\mathbb{R}^d\times\mathbb{R}^d)$, and $I_K:\mathrm{KS}^2(\mathbb{R}^d)\rightarrow \mathrm{KS}^2(\mathbb{R}^d)$ the integral operator induced by $K$. For any bounded subset $B\subset \mathrm{KS}^2(\mathbb{R}^d)$, $I_K(B)$ is relatively compact in $\mathrm{KS}^2(\mathbb{R}^d)$. In fact, consider $\{g_n=I_K(f_n)\}\subset I_K(B)$ and $g\in \mathrm{KS}^2(\mathbb{R}^d)$ such that $g_n\to g$ in $\mathrm{KS}^2(\mathbb{R}^d)$. Since $\{f_n\}$ is bounded in $\mathrm{KS}^2(\mathbb{R}^d)$ there exists $f\in \mathrm{KS}^2(\mathbb{R}^d)$ such that 
$$
\langle f_n, h\rangle_{2}\to \langle f, h\rangle_{2}, \qquad \mbox{as $n\to\infty$, and for all $h\in  \mathrm{KS}^2(\mathbb{R}^d)$}.
$$ 
In particular, for any $x\in\mathbb{R}^d$ it holds that
$$
g_n(x)=\langle f_n, K(x, \,\cdot\,)\rangle_{2}\to \langle f, K(x, \,\cdot\,)\rangle_{2}=I_K(f)(x),\qquad \mbox{as $n\to \infty$},
$$ 
which means that $g=I_K(f)\in I_K(B)$, and $\overline{I_K(B)}$ is compact. 

We essentially argued that if $K\in \mathrm{KS}^2(\mathbb{R}^d\times\mathbb{R}^d)$, then $I_K$ is a compact operator. Furthermore, it is easy to see that 
$$
\langle I_K(f), g\rangle_{KS2}=\langle f, I_K(g)\rangle_{KS2},
$$
for any $f,g\in \mathrm{KS}^2(\mathbb{R}^d)$, then we have established the following.

\begin{cor}
    \label{cor-compact-integral-operator} The integral operator $I_K:\mathrm{KS}^2(\mathbb{R}^d)\rightarrow \mathrm{KS}^2(\mathbb{R}^d)$ generated by $K\in \mathrm{KS}^2(\mathbb{R}^d\times\mathbb{R}^d)$ is a compact self-adjoint operator.
\end{cor}

\section{Mercer's representation of a kernel in the KS-2 space}
\label{Mercer-representation}

A function $K:X\times X\rightarrow\mathbb{R}$, with $X\subset\mathbb{R}^d$ a non-empty set, is called a kernel on $X$. Consider the compact self-adjoint integral operator $I_K:\mathrm{KS}^2(\mathbb{R}^d)\rightarrow \mathrm{KS}^2(\mathbb{R}^d)$ generated a kernel $K\in \mathrm{KS}^2(\mathbb{R}^d\times\mathbb{R}^d)$, given in equation \eqref{integral-operator}. 

The Spectral theorem for compact self-adjoint operators on Hilbert spaces implies that there exists $\{k_n\}\subset\mathrm{KS}^2(\mathbb{R}^d)$ an orthonormal complete system of eigenfunctions of $I_K$ associated with the sequence of positive real numbers $\{\lambda_n\}$ as eigenvalues and such that $\lambda_n\to 0$. For any $f\in \mathrm{KS}^2(\mathbb{R}^d)$, one can write
$$
f=\sum_{n=1}^{\infty}\langle k_n, f\rangle_{KS2} k_n,
$$
in order to obtain
$$
I_K(f)=\sum_{n=1}^{\infty}\lambda_n\langle k_n, f\rangle_{KS2} k_n, \qquad f\in \mathrm{KS}^2(\mathbb{R}^d).
$$
Since $I_K(k_n)(x)=\lambda_nk_n(x)=\langle K(x,\, \cdot \,), k_n\rangle_{KS2}$, for $n=1,2,\ldots$, and 
$$
K(x,\, \cdot \,)=\sum_{n=1}^{\infty}\langle k_n, K(x,\, \cdot \,)\rangle_{KS2} k_n, \qquad x\in \mathbb{R}^d,
$$
we deduce that 
$$
K(x,\, \cdot \,)=\sum_{n=1}^{\infty}\lambda_nk_n(x)k_n, \qquad x\in \mathbb{R}^d.
$$
This delivers us the following Mercer's theorem for kernels $K\in \mathrm{KS}^2(\mathbb{R}^d)$.

\begin{prop} 
    Let $K\in \mathrm{KS}^2(\mathbb{R}^d)$ and $\{k_n\}$ an orthonormal complete system of $\mathrm{KS}^2(\mathbb{R}^d)$ formed by the eigenfunctions associated with the non-negative sequence of eigenvalues $\{\lambda_n\}$ of the integral operator $I_K$ induced by $K$. Then, for any $x\in \mathbb{R}^d$, 
\begin{equation}\label{mercer-representation}
K(x,y)= \sum_{n=1}^{\infty}\lambda_nk_n(x)k_n(y), \qquad   y\in \mathbb{R}^d, \mu-a.e..    \end{equation}
\end{prop} 

It is easy to see that $K$ as in proposition is a definite positive kernel \cite{Berg-2012harmonic}, i.e., besides symmetric it satisfies
$$
\sum_{i,j=1}^nc_ic_jK(x_i, x_j)\geq 0, 
$$
for any $n\in\mathbb{N}$, $\{x_1, \ldots, x_n\}\subset X=\mathbb{R}^d$, and $\{c_1, \ldots, c_n\}\subset\mathbb{R}$. It well known that a definite positive kernel $K$ on $X\neq \emptyset$ is non-negative on the diagonal $\Delta_X:=\{(x,x): x\in X\}$ and it 
satisfies the inequality
\begin{equation}\label{triangular-kernel}
|K(x,y)|\leq K(x,x)K(y,y), \qquad x,y\in  X.    
\end{equation}
Straightly follows that a definite positive kernel on $X$ which is bounded on the diagonal $\Delta_X$ is bounded on $X\times X$.

The Aronszajn theory of Reproducing Kernel Hilbert Spaces (RKHSs) \cite{aron-MR0051437} asserts that for definite positive kernels there exists a unique Hilbert space $(\mathcal{H}_{K},\langle \,\cdot\, , \cdot \, \rangle_K)$ of functions on $\mathbb{R}^d$, satisfying:\\
i. $K(x,\, \cdot\, )\in \mathcal{H}_{K}$ for all $x\in \mathbb{R}^d$;\\
ii. span$\{K(x,\, \cdot\, ) : x\in \mathbb{R}^d\}$  is dense in $\mathcal{H}_{K}$;\\
iii. (Reproducing property) $f(x)=\langle{f},{K(x,\, \cdot\, )}\rangle_{K}$, for all $x\in\mathbb{R}^d$ and $f\in \mathcal{H}_{K}$. \\
We observe that for any $f\in \mathcal{H}_K$, by the reproducing property and the Cauchy-Schwarz inequality, it holds
\begin{eqnarray*}
|f(x)|=\langle f, K(x,\, \cdot\, ) \rangle_K &\leq & \|f\|_K\|K(x,\, \cdot\, )\|_K= \|f\|_K(K(x,x))^{1/2}, \qquad x\in\mathbb{R}^d,
\end{eqnarray*}
and, analogously, for all $x,y\in\mathbb{R}^d$, it holds
\begin{eqnarray*}
    |f(x)-f(y)| \leq  \|f\|_K \|K(x,\, \cdot\, )-K(y,\, \cdot\, )\|_K=  \|f\|_K \left(K(x,x)+K(y,y)-2K(x,y)\right)^{1/2}.
\end{eqnarray*}
Therefore, if $K$ is bounded on the diagonal $\Delta=\{(x,x): x\in \mathbb{R}^d\}$, then the RKHS $\mathcal{H}_K$ contains only bounded functions, and if $K$ is a continuous kernel, then the RKHS $\mathcal{H}_K$ contains only continuous functions.

Due the Mercer's representation \eqref{mercer-representation}, we have that
\begin{equation}\label{RKHSofK}
    \mathcal{H}_{K}=\left\{\sum_{n=1}^\infty c_n\sqrt{\lambda_n}k_n : \{c_2\}\in\ell_2\right\},
\end{equation}
endowed with the inner product given by
\begin{equation}\label{inner-RKHSofK}
    \langle f , g \rangle_K=\sum_{n=1}^{\infty} c_nd_n,
\end{equation}
with $f,g\in \mathcal{H}_{K}$ represented by
\begin{equation*}
f=\sum_{n=1}^\infty c_n\sqrt{\lambda_n}k_n,\qquad g=\sum_{n=1}^\infty d_n\sqrt{\lambda_n}k_n,
\end{equation*}
is the RKHS of the kernel $K\in\mathrm{KS}^2(\mathbb{R}^d)$. This occurs due the unique RKHS associated with the definite positive kernel $K$ satisfying the properties i), ii), and iii) above. In fact, for any $x\in\mathbb{R}^d$ we have that
$$
K(x, \,\cdot\,)= \sum_{n=1}^{\infty}c_n(x)\sqrt{\lambda_n}k_n,
$$
with $d_n(x):=\sqrt{\lambda_n}k_n(x)$, $n=1,2,\ldots$, satisfying 
$$
\sum_{n=1}^{\infty}d_n^2(x)=\sum_{n=1}^{\infty}\lambda_n[k_n(x)]^2=K(x,x),
$$
it means that $\{d_n(x)\}\in\ell_2$, and $K(x, \,\cdot\,)\in \mathcal{H}_{K}$. Also, if $f\in \mathcal{H}_{K}$ is given by 
$$
f=\sum_{n=1}^\infty c_n\sqrt{\lambda_n}k_n,
$$
we have
$$
\langle f, K(\,\cdot \,, x)\rangle_K=\sum_{\lambda_n\neq 0, n=1}^\infty c_nd_n(x)=\sum_{n=1}^\infty c_n\sqrt{\lambda_n}k_n(x)=f(x), \qquad x\in\mathbb{R}^d,
$$
then, the reproducing property holds. 

For any $f\in \mathcal{H}_{K}$ represented by $f=\sum_{n=1}^\infty c_n\sqrt{\lambda_n}k_n$,
the Parseval's identity implies that
$$
\|f\|_{KS2}^2=\sum_{n=1}^{\infty}c_n^2\lambda_n<\infty,
$$
since $\{c_n\}\in\ell_2$ and $\lambda_n\to 0$. Let $\iota_K:\mathcal{H}_K\rightarrow \mathrm{KS}^2(\mathbb{R}^d)$ be the inclusion map. Then, 
$$
\|\iota_K(f)\|_{KS2}^2=\sum_{n=1}^{\infty}c_n^2\lambda_n \leq \sup_n\{\lambda_n\}\sum_{\lambda_n\neq 0, n=1}^{\infty}c_n^2=\sup_n\{\lambda_n\}\|f\|_{K}^2, \qquad f\in \mathcal{H}_K,
$$
and we have established the following.

\begin{prop}\label{embedding}
    Let $K\in \mathrm{KS}^2(\mathbb{R}^d)$ represented as in formula \eqref{mercer-representation}, and $(\mathcal{H}_K, \langle \,\cdot\, , \,\cdot\,\rangle_K)$ the RKHS of $K$ given in \eqref{RKHSofK}. Then, $\mathcal{H}_K\hookrightarrow \mathrm{KS}^2(\mathbb{R}^d)$ compactly. 
\end{prop}

Keeping the notation and under the assumptions of the proposition above, the sequence of eigenvalues $\{\lambda_n\}$ of the integral operator \eqref{integral-operator} defines the multiplier operator $j_K:\mathrm{KS}^2(\mathbb{R}^d)\rightarrow \mathcal{H}_K$ given by 
\begin{equation}\label{j_K}
j_K\left(\sum_{n=1}^{\infty}\langle f, k_n\rangle_{KS2} k_n\right)=\sum_{n=1}^{\infty}\langle f, k_n\rangle_{KS2}\sqrt{\lambda_n}k_n, \qquad f\in \mathrm{KS}^2(\mathbb{R}^d).
\end{equation}
It is clear that $\|j_K\|\leq 1$, and
$$
\|j_K(f)\|_{K} = \|f\|_{KS2}, \qquad f\in  \mathrm{KS}^2(\mathbb{R}^d),  
$$
if and only if $\lambda_n>0$, $n=1,2,\ldots$. 

\subsection{Upper bound for the covering numbers of the embedding $\mathcal{H}_K\hookrightarrow \mathrm{KS}^2(\mathbb{R}^d)$}
\label{covering}

The interest on the embedding $\iota_K$ is to estimate the related covering numbers. The covering numbers of the operator $\iota_K$ are defined in terms of the concept of the covering numbers for metric spaces. If $A$ is a subset of a metric space $X$ and $\epsilon>0$, then the \emph{covering numbers} $\mathcal{C} (\epsilon, A)$ are the minimal number of balls of $X$ with radius $\epsilon$ which covers $A$. We denote by $B_K$ the unit ball in $\mathcal{H}_{K}$ and the \emph{covering numbers of the kernel $K$} (or the covering numbers of the RKHS of the kernel $K$) are defined in terms of the covering numbers of the embedding $\iota_K$, given as follows
$$
\mathcal{C} (\epsilon, \iota_K):= \mathcal{C}(\epsilon,\iota_K(B_K)), \quad \epsilon>0.
$$

We consider $K\in \mathrm{KS}^2(\mathbb{R}^d)$ possessing a Mercer's representation given by 
\begin{equation}\label{exe-mercer-representation}
K(x,y)= \sum_{n=1}^{\infty}\lambda_nY_n(x)Y_n(y), \qquad   x,y\in \mathbb{R}^d,    
\end{equation}
with $\lambda_n>0$, for $n=1,2,\ldots$, satisfying $\lambda_n=O\left(2^{-n(2d+1)}\right)$. We observe that $K$ represented as \eqref{exe-mercer-representation} is a definite positive kernel on $\mathbb{R}^d$ and it holds that
\begin{equation}\label{exe-mercer-representation-ineq}
K(x,x)= \sum_{n=1}^{\infty}\lambda_n(Y_n(x))^2\leq \frac{1}{2\alpha_d^2}\sum_{n=1}^{\infty}\lambda_n2^{n(2d+1)}, \qquad   x\in \mathbb{R}^d,
\end{equation}
since 
$$
Y_k^2=\frac{2^{k-1}}{\lambda^2(B_k)}\chi_{B_k}^2,\qquad \lambda(B_k)= \alpha_d2^{-kd},
$$
for $k=1,2,\ldots$, as defined in Section \ref{sec-KS-space}). The RKHS $(\mathcal{H}_K, \langle \,\cdot\, , \,\cdot\,\rangle_K)$ of $K$ is given by formula \eqref{RKHSofK} replacing the sequence of functions $\{k_n\}$ by $\{Y_n\}$, the orthonomal complete system of $(\mathrm{KS}^2(\mathbb{R}^d), \langle\,\cdot \,, \,\cdot \,\rangle_{KS2})$.  We keep writing $\iota_K$ for the compact embedding $\mathcal{H}_K\hookrightarrow \mathrm{KS}^2(\mathbb{R}^d)$ established in Proposition \ref{embedding}, and we note that \begin{equation}\label{norm-iota}
    \|\iota_K\|\leq \sup_{n}\sqrt{\lambda_n}\leq c2^{-(2d+1)/2},
\end{equation}
if $\lambda_n\leq c^2 2^{-n(2d+1)}$, for $n=1,2\ldots$, and some $c>0$.

The next theorem has its spherical version for continuous and invariant under rotation fixing a pole kernels presented in the Proposition 3.1 of the reference \cite{GONZALEZ2024128121}. We employ the same technique here and it is based on the application of standard properties for covering numbers. The properties stated in a general context and more about the theme can be found in \cite{GONZALEZ2024128121} and references quoted in there. The symbol $\log$ will stand for the logarithm with basis 2.

\begin{thm}\label{thmuppergeo}  Let $K\in \mathrm{KS}^2(\mathbb{R}^d)$ represented as the series expansion \eqref{exe-mercer-representation}, with $\lambda_n>0$, for $n=1,2,\ldots$, satisfying $\lambda_n=O\left(2^{-n(2d+1)}\right)$. Then
$$
\limsup_{\epsilon\to 0^+}\frac{\log (\mathcal{C}(\epsilon, \iota_K ))}{[\log(1/\epsilon)]^{2}}\leq\frac{2}{2d+1} . 
$$
\end{thm}

\pf We consider $V_m=[\sqrt{\lambda_1}k_1, \ldots, \sqrt{\lambda_m}k_m]$, the subspace generated by $\{\sqrt{\lambda_1}k_1, \ldots, \sqrt{\lambda_m}k_m\}$, and $V_m^s$ the orthonormal complement of $V_m$, for $m=1,2,\ldots$. We write $\rho_m$ and $\rho_m^s$ for the orthogonal projections from $\mathcal{H}_k$ onto $V_m$ and $V_m^s$, respectively, for $m=1,2,\ldots$. We also employ the symbols $\rho_m$ and $\rho_m^s$ for the natural embeddings $V_m\hookrightarrow \mathrm{KS}^2(\mathbb{R}^d)$ and $V_m^s\hookrightarrow \mathrm{KS}^2(\mathbb{R}^d)$, respectively.

Let $c>0$ such that $\lambda_n\leq c^22^{-n(2d+1)}$, for $n=1,2\ldots$. As obtained in formula \eqref{norm-iota}, from Proposition \ref{embedding}, we have that
\begin{equation}\label{norm-projs}
    \|\rho_m\|\leq \sup_{n\leq m}\sqrt{\lambda_n}\leq c2^{-(2d+1)/2},\qquad \|\rho_m^s\|\leq \sup_{n>m}\sqrt{\lambda_n}\leq c2^{-(m+1)(2d+1)/2}.
\end{equation}
For $\epsilon >0$, we consider $m:=m(\epsilon)\in \mathbb{Z}_+$ as follows
\begin{equation}\label{alphaMEpsilonexp}
       \frac{c}{2^{(m+1)(2d+1)/2}}< \frac{\epsilon}{2} < \frac{c}{2^{m(2d+1)/2}},
\end{equation} 
by applying the logarithm, we have that 
\begin{equation}\label{m1}
m\approx \frac{2}{2d+1}\log\left(2c/\epsilon\right),  \qquad \mbox{as $\epsilon\to 0^+$}.  
\end{equation}
Due the fact that $\|\rho_m^s\|<\epsilon/2$ by
inequalities \eqref{norm-projs} and \eqref{alphaMEpsilonexp}, a elementary property for the covering numbers guarantees that
$$
\mathcal{C}\left(\epsilon/2 , \rho_m^s \right)=1, \quad m=1,2,\ldots.
$$
Considerations above and the basic property bellow 
$$
\mathcal{C} (\epsilon, \iota_K = \rho_m+\rho_m^s)\leq \mathcal{C} (\epsilon/2, \rho_m)\mathcal{C} (\epsilon/2, \rho_m^s), 
$$
leads us to 
$$
\mathcal{C} (\epsilon, \iota_K)\leq \mathcal{C} (\epsilon/2, \rho_m).
$$
Since $\rho_m$ is an operator on a finite-dimensional Hilbert space it is well known that
$$
\mathcal{C} (\epsilon/2, \rho_m)\leq \left(1+\frac{4\|\rho_m\|}{\epsilon}\right)^{\operatorname{rank}(\rho_m)},
$$
then it follows
$$
\mathcal{C} (\epsilon, \iota_K)\leq \left(1+\frac{4\|\rho_m\|}{\epsilon}\right)^{\operatorname{rank}(\rho_m)}\leq \left(1+\frac{4c}{ 2^{(2d+1)/2}\epsilon}\right)^{m}.
$$
Therefore,
$$
\log(\mathcal{C} (\epsilon, \iota_K))\leq m\log\left(1+\frac{4c}{ 2^{(2d+1)/2}\epsilon}\right)\leq m\log\left(\frac{5c}{ 2^{(2d+1)/2}\epsilon}\right), 
$$
for $m$ large enough. Applying \eqref{m1}, and the asymptotic relation $\log ax\approx \log x$, as $x\to\infty$, for $a>0$, we deduce that
\begin{eqnarray*}
    \log(\mathcal{C} (\epsilon, \iota_K))\leq  m\log\left(\frac{5c}{ 2^{(2d+1)/2}\epsilon}\right)\approx\frac{2}{2d+1}\left(\log\left(1/\epsilon\right)\right)^2, \qquad \mbox{as $\epsilon\to 0^+$},  
\end{eqnarray*}
and it finishes the proof. \eop

\subsection{Lower bound for the covering numbers of the embedding $\mathcal{H}_K\hookrightarrow \mathrm{KS}^2(\mathbb{R}^d)$}\label{lower}

We consider $K\in \mathrm{KS}^2(\mathbb{R}^d)$ possessing a Mercer representation as in \eqref{exe-mercer-representation} with $\lambda_n>0$, for $n=1,2,\ldots$, satisfying $O(\lambda_n)=2^{-n(2d+1)}$. Clearly $K$ represented as \eqref{exe-mercer-representation} is definite positive kernel on $\mathbb{R}^d$ and in order to find a lower estimate we will need the following inequality for the covering numbers of the operator $L:X\rightarrow Y$, with $X$ and $Y$ two n-dimensional Hilbert spaces, 
\begin{equation}\label{CN1}
\sqrt{ \det (L ^* L )}\left(\frac{1}{\epsilon}\right)^{n} \leq \mathcal{C} (\,\epsilon, L ),\quad \epsilon>0.
\end{equation}
You can find the property stated above in \cite{GONZALEZ2024128121}.

\begin{thm}\label{thmlowergeo}  Let $K\in \mathrm{KS}^2(\mathbb{R}^d)$ represented as the series expansion \eqref{exe-mercer-representation}, with $\lambda_n>0$, for $n=1,2,\ldots$, satisfying $O\left(\lambda_n\right)=2^{-n(2d+1)}$. Then
$$
\frac{1}{2 (2d+1)}\leq \liminf_{\epsilon\to 0^+}\frac{\log (\mathcal{C}(\epsilon, \iota_K ))}{[\log(1/\epsilon)]^{2}} . 
$$
\end{thm}

\pf We consider again $V_m = \left[ \sqrt{\lambda_1} k_1 , \ldots, \sqrt{\lambda_m} k_m \right]$, the subspace generated by $\{ \sqrt{\lambda_1} k_1 , \ldots, \sqrt{\lambda_m} k_m \}$ for $m=1,2,\ldots,$ and the composition operator
$$
T_m : V_m \stackrel{i_m}{\longrightarrow} \mathcal{H}_K \stackrel{\iota_K}{\longrightarrow} \mathrm{KS}^2 (\mathbb{R}^d) \stackrel{j_K}{\longrightarrow} \mathcal{H}_K \stackrel{P_m}{\longrightarrow} F_m,
$$
for $m=1,\ldots,$ where $i_m$ stands for the embedding given by $V_m \hookrightarrow\mathcal{H}_K $, $ j_K$ as in \eqref{j_K}, and $P_m$ is the orthogonal projection of $\mathcal{H}_K$ on $F_m := j_K \iota_K i_m (V_m)$. \\
Thereby
$$ 
\langle T_m^* T_m \sqrt{\lambda_r} k_r, \sqrt{\lambda_l} k_l \rangle_K =\langle T_m \sqrt{\lambda_r} k_r , T_m \sqrt{\lambda_l} k_l\rangle_K = \langle \lambda_r k_r, \lambda_l k_l \rangle_K = \lambda_r \delta_{rl}, 
$$ 
and the representing matrix of $T_m^* T_m$ is the following diagonal matrix

\begin{equation}
\begin{pmatrix}
\lambda_1 & 0 & \cdots & 0\\
0 & \lambda_2 & \cdots & 0\\
\vdots & \vdots & \ddots & \vdots\\
0 & 0 & \cdots & \lambda_m
\end{pmatrix}
\end{equation}
for $m=1,2,\ldots,$ whence 
\begin{equation}\label{prod}
\det (T_m^* T_m)= \prod_{k=1}^{m} \lambda_k .
\end{equation}
Due to that $\Vert P_m\Vert=\Vert i_m\Vert=1$ and $\Vert j_K\Vert \leq 1$, we know $\mathcal{C}\left(1 , P_m \right)=\mathcal{C}\left(1 , j_K \right)=\mathcal{C}\left(1 , i_m \right)=1$, and by covering numbers property, we have that for any $\epsilon>0$,

\begin{eqnarray*}
\mathcal{C} (\,\epsilon, T_m ) & = &\mathcal{C} (\,\epsilon, P_m \,j_K \,\iota_K \,i_m ) \\ & \leq & \mathcal{C} (1, P_m ) \mathcal{C} (1,j_K ) \mathcal{C} (\epsilon, \iota_K )\mathcal{C} (1, i_ m )\\ & = & \mathcal{C} (\epsilon, I_K ), \quad m=1,2,\ldots.
\end{eqnarray*}
The lower bound for the covering number in inequality \eqref{CN1}, leads us to
$$
\sqrt{ \det (T_m ^* T_m )}\left(\frac{1}{{}\,\epsilon}\right)^{m} \leq \mathcal{C} (\,\epsilon, T_m  ) \leq \mathcal{C} (\epsilon, I_K ), 
$$
and due to the formula \eqref{prod}, we obtain
$$
\prod_{k=1}^{m} \lambda_k \left(\frac{1}{{}\,\epsilon}\right)^{m} \leq \mathcal{C} (\epsilon, I_K ), \quad \epsilon>0, \quad m=1,2,\ldots.
$$
Applying the natural logarithm in both sides of inequality above, we conclude that for any $\epsilon>0$, 
\begin{equation}\label{lower1geo}
     \frac{1}{2}\sum_{k=0}^{m} \log \left(\lambda_k \right) - m\,\log(\epsilon)\leq \log(\mathcal{C}(\epsilon, I_K )), \quad m=1,2,\ldots.
\end{equation}
For $a>0$ and $m=1,2,\ldots, $ the estimate
$$
\frac{1}{a\, 2^{m(2d+1)}} \leq \lambda_m \leq \lambda_k, \quad k=1,2,\ldots,m.
$$ 
implies 
\begin{equation}\label{ef}
- m^2 \left(\frac{2d+1}{2}\right)+m \log\left(\frac{1}{\sqrt{a}\,\epsilon}\right) \leq \log(\mathcal{C}(\epsilon, I_K )).
\end{equation}
The final steps of the proof are to solve the optimization problem given by the following quantity, 
$$
\varphi_\epsilon (m):= - m^2 \left(\frac{2d+1}{2}\right)+m \log\left(\frac{1}{\sqrt{a}\,\epsilon}\right) 
$$
with $\epsilon>0$ and $m=1,2,\ldots$. 
The critical point of $\varphi_\epsilon$, as a function of $\mathbb{R}$, is given by
$$
c_0 = -\frac{\log(\sqrt{a}\,\epsilon)}{2d+1}.
$$
Since $\lfloor c_0 \rfloor \approx  c_0 \approx \lceil c_0 \rceil $, we can estimates \eqref{ef} in terms of $\varphi_{\epsilon}(c_0)$, even if $c_0$ is not an integer. A simple calculation leads us to 
$$
\varphi_\epsilon (c_0)= \frac{1}{2(2d+1)} (\log \left(1/\sqrt{a}\,\epsilon)\right)^2 \approx  \frac{1}{2(2d+1)} (\log \left(1/\epsilon)\right)^2 , 
$$
and we conclude that 
$$
\frac{1}{2(2d+1)} (\log \left(1/\epsilon)\right)^2 \leq \log(\mathcal{C}(\epsilon, I_K )) ,
$$
for any $\epsilon>0$ and the proof follow.
\eop

\begin{cor}\label{4.5}
    Let $K\in KS^2 (\mathbb{R}^d)$ represented as the series expansion \eqref{exe-mercer-representation}, with $\lambda_n > 0$, for $n=1,2,\ldots,$ satisfying  $\lambda_n \asymp 2^{-n(2d+1)} $, then
\[
\log(\mathcal{C}(\epsilon,I_K ))\asymp {\log\left(\frac{1}{\epsilon}\right)}^{2} ,
\]
as $\epsilon\rightarrow 0^+$.
\end{cor}

Estimates of the asymptotic behavior of the upper bound and lower bound for the covering numbers of unite balls of RKHS have also been recently investigated in other contexts in \cite{GJ,GONZALEZ2024128121}. There, the authors used the RKHS generated by the Weierstrass fractal kernel on a closed interval $I$, and two general classes of RKHS on the unit sphere $\mathbb{S}^d$ of $\mathbb{R}^{d+1}$, where the embeddings focus were $\mathcal{H}_K \hookrightarrow C(I)$ and $\mathcal{H}_K \hookrightarrow C(\mathbb{S}^d)$, respectively, with $C(I)$ and $C(\mathbb{S}^d)$ the space of continuous functions. Versions of Corollary \ref{4.5} can also be found in \cite{GJ,GONZALEZ2024128121}.\\

\noindent\textbf{Acknowledgements.} his research was supported by the São Paulo Research Foundation (FAPESP) under grant numbers 2022/11032-0 and  2021/12213-5. Additional support was provided by the Morá Miriam Rozen Gerber fellowship from the Weizmann Institute of Science. A preliminary version of this work was presented in abstract form at Analysis and PDE in Latin America.

\bibliographystyle{plain}
\bibliography{KSrefsNOVO.bib}

\begin{thebibliography}{10}

\bibitem{converse-regular}
F.~Andrade~da Silva, M.~Federson, R.~Grau, and E.~Toon.
\newblock Converse {L}yapunov theorems for measure functional differential
  equations.
\newblock {\em J. Differ. Equ.}, 286:1--46, 2021.

\bibitem{converse-uniform}
F.~Andrade~da Silva, M.~Federson, and E.~Toon.
\newblock Stability, boundedness and controllability of solutions of measure
  functional differential equations.
\newblock {\em J. Differ. Equ.}, 307:160--210, 2022.

\bibitem{aron-MR0051437}
N.~Aronszajn.
\newblock Theory of reproducing kernels.
\newblock {\em Trans. Amer. Math. Soc.}, 68:337--404, 1950.

\bibitem{Modern}
Robert~G. Bartle.
\newblock {\em A modern theory of integration}, volume~32 of {\em Graduate
  Studies in Mathematics}.
\newblock Amer. Math. Soc., Providence, RI, 2001.

\bibitem{bongiorno}
Donatella Bongiorno and Giuseppa Corrao.
\newblock An integral on a complete metric measure space.
\newblock {\em Real Anal. Exchange}, 40(1):157--178, 2014/15.

\bibitem{B-F-M}
E.~M. Bonotto, M.~Federson, and J.~G Mesquita.
\newblock {\em Generalized ordinary differential equations in abstract spaces
  and applications}.
\newblock John Wiley \& Sons, Hoboken, NJ, 2021.

\bibitem{smale-MR1864085}
Felipe Cucker and Steve Smale.
\newblock On the mathematical foundations of learning.
\newblock {\em Bull. Amer. Math. Soc. (N.S.)}, 39(1):1--49, 2002.

\bibitem{FollandReal}
Gerald~B. Folland.
\newblock {\em Real analysis}.
\newblock Pure and Applied Mathematics (New York). John Wiley \& Sons, Inc.,
  New York, second edition, 1999.
\newblock Modern techniques and their applications, A Wiley-Interscience
  Publication.

\bibitem{gill2009}
Tepper~L. Gill and Woodford~W. Zachary.
\newblock Banach spaces for the {F}eynman integral.
\newblock {\em Real Anal. Exchange}, 34(2):267--310, 2009.

\bibitem{gill2016functional}
Tepper~L. Gill and Woodford~W. Zachary.
\newblock {\em Functional analysis and the {F}eynman operator calculus}.
\newblock Springer, Cham, 2016.

\bibitem{GONZALEZ2024128121}
Karina Gonzalez and Thaís Jordão.
\newblock A close look at the entropy numbers of the unit ball of the
  reproducing hilbert space of isotropic positive definite kernels.
\newblock {\em Journal of Mathematical Analysis and Applications},
  534(2):128121, 2024.

\bibitem{Gordon}
Russell~A. Gordon.
\newblock {\em The integrals of {L}ebesgue, {D}enjoy, {P}erron, and
  {H}enstock}, volume~4 of {\em Graduate Studies in Mathematics}.
\newblock Amer. Math. Soc., Providence, RI, 1994.

\bibitem{Henstock-First}
Ralph Henstock.
\newblock {\em Theory of integration}.
\newblock Butterworths, London, 1963.

\bibitem{kolmogorov-MR0112032}
A.~N. Kolmogorov and V.~M. Tihomirov.
\newblock {$\varepsilon $}-entropy and {$\varepsilon $}-capacity of sets in
  function spaces.
\newblock {\em Uspehi Mat. Nauk}, 14(2(86)):3--86, 1959.

\bibitem{theories}
Douglas~S. Kurtz and Charles~W. Swartz.
\newblock {\em Theories of integration}, volume~13 of {\em Series in Real
  Analysis}.
\newblock World Scientific Publishing Co. Pte. Ltd., Hackensack, NJ, second
  edition, 2012.
\newblock The integrals of Riemann, Lebesgue, Henstock-Kurzweil, and McShane.

\bibitem{Kurzweil-First}
Jaroslav Kurzweil.
\newblock Generalized ordinary differential equations and continuous dependence
  on a parameter.
\newblock {\em Czechoslovak Math. J.}, 7(82):418--449, 1957.

\bibitem{Li-MR1733160}
Wenbo~V. Li and Werner Linde.
\newblock Approximation, metric entropy and small ball estimates for {G}aussian
  measures.
\newblock {\em Ann. Probab.}, 27(3):1556--1578, 1999.

\bibitem{Minh2010-MR2677883}
Ha~Quang Minh.
\newblock Some properties of {G}aussian reproducing kernel {H}ilbert spaces and
  their implications for function approximation and learning theory.
\newblock {\em Constr. Approx.}, 32(2):307--338, 2010.

\bibitem{reid1975anatomy}
WT~Reid.
\newblock Anatomy of the ordinary differential equation.
\newblock {\em The American Mathematical Monthly}, 82:971--984, 1975.

\bibitem{Berg-2012harmonic}
C.~van~den Berg, J.P.R. Christensen, and P.~Ressel.
\newblock {\em Harmonic Analysis on Semigroups: Theory of Positive Definite and
  Related Functions}.
\newblock Graduate Texts in Mathematics. Springer, New York, 2012.

\bibitem{WSS-MR1873936}
Robert~C. Williamson, Alex~J. Smola, and Bernhard Sch\"{o}lkopf.
\newblock Generalization performance of regularization networks and support
  vector machines via entropy numbers of compact operators.
\newblock {\em IEEE Trans. Inform. Theory}, 47(6):2516--2532, 2001.

\bibitem{Zhou-2003-MR1985575}
Ding-Xuan Zhou.
\newblock Capacity of reproducing kernel spaces in learning theory.
\newblock {\em IEEE Trans. Inform. Theory}, 49(7):1743--1752, 2003.

\end{thebibliography}
\bigskip

\medskip
{\small 

\noindent \textsc{F. Andrade da Silva} \\
Departament of Mathematics\\
Universidade de S\~{a}o Paulo\\
13566-590 S\~{a}o Paulo, Brazil\\
\textsc{Email}: ffeandrade@usp.br

\medskip
\noindent {\textsc{K. Gonzalez}}\\
Faculty of Mathematics and Computer Science\\
Weizmann Institute of Science\\
7632706 Rehovot, Israel\\
\textsc{Email}: karina.navarro-gonzalez@weizmann.ac.il

\medskip

\noindent \textsc{T. Jord\~{a}o}\\
Departament of Mathematics\\
Universidade de S\~{a}o Paulo\\
13566-590 S\~{a}o Paulo, Brazil\\
\textsc{Email}: tjordao@icmc.usp.br

\bigskip

\noindent{\textsc{Keywords}:} isotropic kernel; covering numbers; Reproducing Kernel Hilbert Space

\medskip

\noindent{\textsc{MSC}:} 47B06; 42C10; 46E22; 41A60

\end{document}